\documentclass{amsart}

\usepackage[colorlinks=true, urlcolor=blue,bookmarks=true,bookmarksopen=true,citecolor=blue,hypertex]{hyperref}
\usepackage{graphicx}
\usepackage{enumerate}
\usepackage{parskip}

\usepackage{amscd}
\usepackage{amssymb}
\usepackage{amsxtra}
\usepackage{amsmath}
\usepackage{enumerate}
\usepackage{mathrsfs}

\usepackage{color}

\usepackage{amsmath,amsfonts,amssymb}
\usepackage{verbatim}
\usepackage{epsfig}
\usepackage{amsthm}
\usepackage{bbold}

\usepackage[T1]{fontenc}

\theoremstyle{plain}

\newtheorem{thm}{Theorem}

\newtheorem{prop}[thm]{Proposition}

\theoremstyle{definition}
\newtheorem{dfn}[thm]{Definition}

\theoremstyle{remark}

\newtheorem{rmk}[thm]{Remark}

\newtheorem*{qsnonum}{Question}

\theoremstyle{plain}

\newcommand{\Qed}{\hfill \qedsymbol \medskip}

\newcommand{\Id}{{{\mathchoice {\rm 1\mskip-4mu l} {\rm 1\mskip-4mu l}
      {\rm 1\mskip-4.5mu l} {\rm 1\mskip-5mu l}}}}
\newcommand\dist{\operatorname{dist}}

\def\ZZ{\mathbb{Z}}

\def\RR{\mathbb{R}}
\def\CC{\mathbb{C}}

\def\Id{\mathbb{1}}

\def\dd{\mathrm{d}}
\DeclareMathOperator{\Ham}{Ham}

\DeclareMathOperator{\Diff}{Diff}

\DeclareMathOperator{\Fix}{Fix}
\DeclareMathOperator{\Spec}{Spec}

\DeclareMathOperator{\Area}{Area}
\DeclareMathOperator{\Arg}{Arg}
\DeclareMathOperator{\w}{wind}

\definecolor{lgray}{rgb}{0.5,0.5,0.5}

\begin{document}

\pagestyle{headings}

\bibliographystyle{alphanum}

\title[Hofer gap between integrable and autonomous Hamiltonians]{A gap in the Hofer metric between integrable and autonomous Hamiltonian diffeomorphisms on surfaces}

\date{\today} 

\author{Michael Khanevsky}
\address{Michael Khanevsky, Faculty of Mathematics,
Technion - Israel Institute of Technology
Haifa, Israel}
\email{khanev@technion.ac.il}

\begin{abstract}
Let $\Sigma$ be a compact surface equipped with an area form. There is an long standing open question by Katok,
which, in particular, asks whether every entropy-zero Hamiltonian diffeomorphism of a surface 
lies in the $C^0$-closure of the set of integrable diffeomorphisms.
A natural generalization of this question is to ask to what extent one family of `simple' Hamiltonian diffeomorphisms of $\Sigma$ can be approximated by the other.
In this paper we show that the set of autonomous Hamiltonian diffeomorphisms is not Hofer-dense
in the set of integrable Hamiltonians. We construct explicit examples of integrable diffeomorphisms which cannot be Hofer-approximated
by autonomous ones.
\end{abstract}

\maketitle

\section{Introduction}

There is the following long standing open question of Katok:
"\emph{In low dimensions, is every volume-preserving dynamical system with zero topological entropy a limit of integrable systems?}"
This is stated as Problem 1 in ~\cite{Katok}, but relates also to the much older paper ~\cite{An-Kat}. 
Low dimensions means maps on surfaces or flows on $3$-dimensional manifolds. 
 
A straightforward generalization of this question is to ask to what extent one `simple' family of Hamiltonian diffeomorphisms 
can be approximated by the other. While the reader is welcome to substitute his own definitions for a family of \emph{simple} diffeomorphisms
and \emph{approximation}, some natural choices would be to consider autonomous/integrable/entropy-zero maps and an approximation in terms of a
commonly used metric on the Hamiltonian group.
While the original question of Katok is still open, some of these generalizations are easier to access and still yield an interesting and non-obvious statement.
In this paper we show that the set of autonomous Hamiltonian diffeomorphisms on a symplectic surface is not Hofer-dense
in the set of integrable Hamiltonians. We construct explicit examples of integrable diffeomorphisms which cannot be Hofer-approximated
by autonomous ones.

\medskip

Let $\Sigma$ be a compact surface with an area form $\omega$ and
denote by $\Ham(\Sigma)$ the group of Hamiltonian diffeomorphisms of $(\Sigma, \omega)$ that are compactly supported in the interior of $\Sigma$. For every smooth function
$H\colon\Sigma\to\RR$ compactly supported in the interior, there exists a unique vector field $X_H$ which satisfies
$$
dH(\cdot)=\omega(X_H,\cdot).
$$
It is easy to see that $X_H$ is tangent to the level sets of $H$. Let $h$ be the
time-one map of the flow $h_t$ generated by $X_H$. The diffeomorphism $h$ is
area-preserving, it belongs to $\Ham(\Sigma)$ and every Hamiltonian diffeomorphism arising this way is called
{\em autonomous}. Such a diffeomorphism is easy to understand
in terms of its generating Hamiltonian function.  More generally, autonomous maps $h \in \Diff (\Sigma)$ appear as time-$t$ maps
of time-independent flows and also have relatively simple dynamics.

A diffeomorphism $g \in \Diff (\Sigma)$ is \emph{integrable} if it admits a smooth invariant function $H : \Sigma \to \RR$
(a \emph{first integral}) such that $H$ is not constant on any open set.
Such a function $H$ induces a lamination on $\Sigma$ by its level sets. This lamination is regular away from the critical points of $H$ and 
is invariant under $g$.

If $g$ is an autonomous Hamiltonian generated by $H : \Sigma \to \RR$, the function $H$ serves as an integral for $g$. If $H$ is constant on certain
open sets (for example, outside its support) one may perturb $H$ there. Such a perturbed function is still invariant under $g$ since the perturbation occurs 
in the set of fixed points of $g$. Not all integrable Hamiltonians are autonomous, an example is the diffeomorphism $g$ constructed in the proof of
Theorem~\ref{T:main}.
Note that a non area-preserving autonomous diffeomorphism $g \in \Diff (\Sigma)$ might be non-integrable 
(e.g. a flow with an attracting point), yet it admits an invariant lamination by flow lines.

While we restrict our attention to Hamiltonian diffeomorphisms on surfaces, there is a hierarchy
\[
  \Ham_{Aut}(\Sigma) \subset \Ham_{Int} (\Sigma) \subset \Ham_{Ent} (\Sigma)
\]
of autonomous / integrable / entropy-zero Hamiltonians on $\Sigma$. 

We consider the following version of Katok's question:
\begin{qsnonum}
Does every integrable Hamiltonian on $\Sigma$ lie in the Hofer-closure of the set of autonomous Hamiltonian diffeomorphisms?
\end{qsnonum}

In this paper we answer in negative the above question for every compact surface $\Sigma$. Our main result is the following

\begin{thm} \label{T:main}
There exists an integrable $g \in \Ham(\Sigma)$ which cannot be presented as a Hofer-limit of autonomous Hamiltonian diffeomorphisms.
\end{thm}

In Section~\ref{S:Braids} we describe a class of braids in the unit disk $D$ which cannot be induced by periodic trajectories of an autonomous flow.
Section~\ref{S:Proof} constructs $g \in \Ham_{Int} (D)$ by composition of two autonomous Hamiltonian diffeomorphisms. We show that
$g$ induces a `non-autonomous' braid $B$ by its periodic trajectories.
In Section~\ref{S:Hofer} we discuss an adaptation of the braid stability result of \cite{Al-Me:braid}. It implies that Hofer-small perturbations
of $g$ that are non-degenerate also induce $B$ by periodic trajectories. Theorem~\ref{T:main} follows for the disk $D$ by a perturbation argument.
Section~\ref{S:Surf} adapts this argument to general surfaces.

\medskip

Related results: this paper is a follow-up for ~\cite{Bra-Kha:C0-gap} which provides a negative answer for the similar question where the Hofer metric 
is replaced by the $C^0$ one. Moreover, in the case of the annulus $\Sigma = S^1 \times [0, 1]$ Theorem~2 of \cite{Bra-Kha:C0-gap} shows
\begin{thm} \label{T:annulus}
Given $R>0$, there exists $g_R \in \Ham_{Int}(S^1 \times [0, 1])$ such that the Hofer distance $\dist_H (g_R, \Ham_{Aut}) > R$.
\end{thm}
While this result is much stronger than Theorem~\ref{T:main}, the tools used in the proof (the Calabi quasimorphisms of \cite{En-Po:calqm}) 
are limited to the case of the annulus and cannot be applied to a general $\Sigma$.

\cite{Po-Sh:egg-beater} constructs, among other results, chaotic (positive-entropy) Hamiltonian diffeomorphisms on higher-genus surfaces
which are Hofer-far from $\Ham_{Aut}$. More examples are provided in \cite{Ch-Me:Hof-entr}, \cite{Ch:eggbeater}.
\cite{Al-Me:braid} shows that $\Ham_{Ent}(\Sigma)$ is Hofer-closed in $\Ham(\Sigma)$.

This paper follows the lines of \cite{Al-Me:braid}, using its methods and tools.
We extend the braid stability results to certain degenerate Hamiltonian diffeomorphisms (this is straightforward) and describe examples 
of braids that are induced by integrable but not by autonomous maps, which is the main contribution of this text.

\subsection*{Acknowledgments.}
We thank Marcelo Alves for explanation of the braid stability results.

\section{Braids and autonomous flows}\label{S:Braids}

Denote by $B_k$ the pure braid group with $k$ strands in the the unit disk $D$.
Given a flow $g_t \in \Diff (D), t \in [0,1]$, we consider the set $\Fix(g_1)$ of fixed points of $g_1$
(these are $1$-periodic points under iterations of $g_t$). Trajectories under $g_t$ of a $k$-tuple of distinct points 
$\mathbf{p} = (p_1, \ldots, p_k)$ from $\Fix(g_1)$ define an element $B_\mathbf{p} \in B_k$.
(Strictly speaking, this defines only a conjugation equivalence class in $B_k$. But in order to simplify the statements we will
still call it a `braid', keeping in mind ambiguity up to a conjugation).

Denote by $B_k^{Aut} \subset B_k$ the set of \emph{autonomous braids} that are induced by autonomous flows $g_t$.
We analyze the structure of braids in $B_k^{Aut}$, which helps to define an obstruction for a braid to be autonomous.
As $\Ham(D)$ is contractible, $g \in \Ham(D)$ determines the homotopy class $\widetilde{g} \in \widetilde{\Ham(D)}$. Hence a choice of a $k$-tuple
$\mathbf{p} \subset \Fix(g)$ induces (up to a conjugation) a braid $B_\mathbf{p} \in B_k$. In the case $B_\mathbf{p} \notin B_k^{Aut}$ 
this serves as a proof that $g$ is not autonomous.

Given a flow $g_t$ and two distinct  points $p, q \in \Fix(g_1)$, consider the \emph{winding number} (a.k.a. the linking number) 
of their parametrized trajectories:
\[
  \w_g (p, q) = \deg \left(\Arg (g_t(p) - g_t (q))\right)
\]
where $\Arg (g_t(p) - g_t (q))$ is seen as a function $[0, 1] / \{ 0, 1\} \to S^1$.
The winding number is symmetric:
\[
  \w_g (p, q) = \w_g (q, p)
\]
and additive: given $p, q \in \Fix(g_1) \cap \Fix(h_1)$, the winding number for the composition of flows
\[
  \w_{g \# h} (p, q) = \w_g (p, q) + \w_h (p, q).
\]

We note that winding numbers descend to braids $B_\mathbf{p}$ induced by periodic trajectories of the flow:
denote by $B_{p,q}$ the sub-braid of $B_\mathbf{p}$ formed by the strings of $p$ and $q$.
Then $\w_g (p, q)$ corresponds to the natural isomorphism $B_2 \to \ZZ$ applied to $B_{p,q}$.
The winding numbers are preserved under conjugations, hence are not affected by ambiguity of 
the definition of $B_\mathbf{p} \in B_k$.

Given an $h \in \Ham(D)$, it determines the lift $\widetilde{h} \in \widetilde{\Ham(D)}$ to the universal cover, 
hence allows definition of winding numbers for points in $\Fix(h)$ which does not depend on the isotopy in $\Ham(D)$ connecting
$h$ to the $\Id_D$.

\medskip

Suppose $g_t$ is an autonomous flow in $\Diff(D)$. 
$1$-periodic points of autonomous flows arise in two possible scenarios: as an equilibrium point (they are fixed by $g_t$ for all $t \in \RR$) or 
a point $p$ whose trajectory $c_p$ is a simple loop in $D$. Moreover, in this case all points on $c_p$ are $1$-periodic as well.
(Here the notion \emph{trajectory} refers to the image of the curve $\{g_t(p)\}_t$ in $D$, forgetting the parametrization.)
The trajectory $c_p$ of a non-equilibrium point $p$ bounds a smooth open disk $D_{c_p}$ in $D$. 
In order to keep the notation uniform, we will think of constant trajectories $c_q$ of an equilibrium point $q$ that they bound the singleton set $\{q\}$ 
(disk of radius zero centered at $q$).

We introduce a partial order on the set of $1$-periodic trajectories of $g_t$:
$c \leq c'$ if $D_c \subset D_{c'}$. In other words, $c'$ either encircles $c$ or coincides with it. It is easy to verify that this 
defines a weak partial order. ($c$, $c'$ are incomparable if and only if they are `separated' -- neither one encircles the other.)

One can see that $\w_g (p, q) = 0$ if and only if the trajectories $c_p, c_q$ are incomparable. Indeed, two separated trajectories can be 
homotoped to constant ones preserving their winding numbers. If $c_q \lneqq c_p$, then $p$ is not an equilibrium point and $c_q$ can be contracted 
inside the disk $D_{c_p}$. As $p$ is not an equilibrium, it must traverse $c_p$ at least once. Hence the path of $p$ encircles $c_q$ (and also its contracted version) 
at least once.
If $c_q = c_p$ (namely, $p, q$ share the same trajectory under $g_t$), $\w_g (p, q)$ agrees with the winding number of $\{g_t (p)\}$
around interior points of $D_{c_p}$ and cannot be zero. (Once again, it can be seen by verifying that $p \neq q \in c_p$ implies they are not equilibrium 
and contracting the trajectory of $q$ inside $D_{c_p})$.
If $p=q$ the winding number $\w_g (p, q)$ is not defined.

Moreover, a similar argument implies that if three distinct $1$-periodic points $p, q, q'$ satisfy $c_{q} \leq c_p$ and $c_{q'} \leq c_p$ then $\w_g (p, q) = \w (p, q')$.

\medskip

\begin{figure}[!htbp]  
\begin{center}
\includegraphics[width=0.5\textwidth]{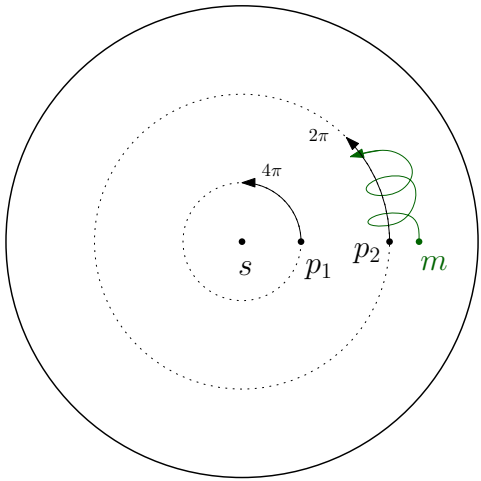}
\caption{The non-autonomous braid $B_{na}$.}   \label{F:motion}
\end{center}
\end{figure}

\textbf{Non-autonomous braids:} we present an example $B_{na}$ of a braid 
(a `non-auto\-nomous' braid) which cannot appear under an autonomous flow in $D$.
This parti\-cular example easily extends to a large class of braids.
Consider the motion as in the Figure~\ref{F:motion} (`a sun, two planets and a moon')
of four $1$-periodic points $s$ (the `sun'), $p_1, p_2$ (`planets') and $m$ (the `moon' of $p_2$).
In time $1$ the planet $p_1$ does two turns around the sun, $p_2$ one turn and the moon three turns around $p_2$ following the orbit of $p_2$ around $s$.
One can easily describe a flow $g_t \in \Diff (D)$ which induces the braid $B_{na}$ defined by this motion. Moreover, 
in Section~\ref{S:Proof} we explicitly construct $g_t$ as a Hamiltonian flow whose time-$1$ map is integrable.
The winding numbers of the strands satisfy
\[
 \begin{split}
	\w (p_1, s) &= 2,  \\
	\w (p_2, s) &= \w (p_2, p_1) = 1, \\
	\w (m, s) &= \w (m, p_1) = 1, \quad \w (m, p_2) = 3 .
 \end{split}
\]
We claim that this choice of the winding numbers obstructs existence of a partial order of trajectories, hence $g_t$ is not autonomous. 
Assume by contradiction that this is not the case and trajectories admit a partial order.
As the winding numbers are non-zero, any two trajectories from $\{c_s, c_{p_1}, c_{p_2}, c_m\}$ are comparable under the partial order. 
(Here we must keep in mind the possibility that some trajectories may coincide.)
This implies that the set $\{c_s, c_{p_1}, c_{p_2}, c_m\}$ is linearly ordered. 
Therefore it has a maximal element and a periodic point which belongs to this `maximal' trajectory has the same winding number
with each of the remaining three periodic points. However, neither point satisfies this condition, which leads to a contradiction.

\begin{rmk}
  This example easily extends to a more general case. First and most obvious, any braid $B$ which contains $B_{na}$ as a subbraid
  by `forgetting' some of the strands of $B$, is non-autonomous.
  
  Looking in the argument deeper, we obtain the following simple obstruction for a braid $B$ to be autonomous: 
  enough to find a subbraid $B'$ with $k$ strands such that 
  \begin{enumerate}
	  \item 
		 the winding number of any two strands in $B'$ is non-zero,
	  \item
		 none of the strands of $B'$ is `maximal': for strands $s_1, s_2$ of $B'$ the winding number $\w (s_1, s_2)$ is not constant in 
		 $s_2$ for any choice of $s_1$.
  \end{enumerate}
  
  Non-autonomous braids which satisfy this obstruction exist already in $B_3$ (consider a braid with three strands and winding numbers $\{1, 2, 3\}$), 
  but $B_{na}$ is the simplest example that the author was able to induce by an \emph{integrable} motion.
\end{rmk}

\section{Hofer's metric and stability of braids}\label{S:Hofer}

Let $(M, \omega)$ be a symplectic manifold, $g$ a Hamiltonian diffeomorphism with compact support in $M$.
The Hofer norm $\|g\|$ (see ~\cite{Hof:TopProp}) is defined by 
\[
  \|g\| = \inf_G \int_0^1 \max \left(G(\cdot, t) - \min G(\cdot, t) \right) \dd t 
\]
where the infimum goes over all compactly supported Hamiltonian functions $G : M \times [0,1] \to \RR$
such that $g$ is the time-1 map of the induced flow. The Hofer distance is given by
\[
  \dist_H (g_1, g_2) = \| g_1 g_2^{-1}\|.
\]
and, intuitively, it expresses the amount of mechanical energy needed to deform the dynamics of $g_1$ into $g_2$.

\medskip

In ~\cite{Al-Me:braid} the authors show the following: suppose a non-degenerate $f \in \Ham (D)$ 
is generated by a flow $\phi_t \in \Ham (D)$. (Non-degenerate for $f \in \Ham(D)$ means $f$ is non-degenerate in 
the interior of the support of $f$.)
Let $B = (\gamma_1, \ldots, \gamma_k)$ be a braid of distinct $1$-periodic orbits of $\phi_t$. 
Then under certain regularity assumptions on $f$, any Hofer-small and non-degenerate perturbation $f'$ of $f$ will form a copy
of the braid $B$ by trajectories of its fixed points.
Very roughly, if the energy difference $\dist_H (f, f')$ is too small to disrupt any individual orbit $\gamma_i$ 
(that is, the perturbation cannot annihilate $\gamma_i$, and, in the case of $\Sigma = S^2$ instead of $D$, no bubbling can occur) 
then the perturbation cannot disrupt the braid $B$ as a whole.

More precisely, the paper suggests the following assumptions on $f \in \Ham (D)$ which is induced by a compactly-supported 
Hamiltonian function $F_t : D \to \RR$ and a braid of periodic points $B = (\gamma_1, \ldots, \gamma_k)$:
\begin{dfn} \label{D:isolated}
  (Definition~2.1 in ~\cite{Al-Me:braid}) suppose $f$ is non-degenerate and let $\epsilon > 0$. 
  $B$ is \emph{$\epsilon$-isolated} for the Hamiltonian action $\mathcal{A}_F$ if 
  all $a, b \in \Spec (f) = \mathcal{A}_F (\Fix (f))$ satisfy either $a = b$ or $|a - b| \geq \epsilon$. 
  We stress that distinct periodic orbits are allowed to have the same action, but if the actions are different, they differ at least by $\epsilon$.
  
  In addition there is a rather technical assumption that if $a \in \mathcal{A}_F (B)$, then all orbits of $\Fix(f)$ with action $a$ appear in $B$.
\end{dfn}
In the case of a closed surface $\Sigma$, ~\cite{Al-Me:braid} proposes a less restrictive
\begin{dfn} \label{D:quasi-isolated} (Definition~3.3 in the text)
  suppose $f \in \Ham(\Sigma)$ is a non-degenerate diffeomorphism generated by a Hamiltonian function $F_t : \Sigma \to \RR$.
  Pick a compatible almost complex structure $J$ and $\epsilon > 0$. Then $B$ is \emph{$\epsilon$-quasi-isolated} for 
  $\mathcal{A}_F$ if 
  \begin{itemize}
	  \item
		  $\forall \gamma, \gamma' \in B$, $|\mathcal{A}_F (\gamma) - \mathcal{A}_F (\gamma')|$ is either $0$ or $\geq \epsilon$.
	  \item
		  Floer trajectories converging on either end to on orbit from $B$ are either constant or have energy $\geq \epsilon$.
  \end{itemize}
(For a detailed explanation and recollections from the Floer theory we refer the reader to the text \cite{Al-Me:braid}.)
\end{dfn}

In the proof of the braid stability the assumptions of Definition~\ref{D:isolated} translate to the following.
Select appropriate auxiliary Floer data and consider filtered Floer complexes $CF^{(\kappa - \delta, \kappa+\delta)}$
for some $\kappa \in \mathcal{A}_F (B)$ and $\delta > 0$. Then (Step~1 in the proof of Theorem~3)
\begin{itemize}
  \item 
	  $CF^{\left(\kappa - \frac{\epsilon}{100}, \kappa+\frac{\epsilon}{100}\right)} (F) \simeq 
	  CF^{\left(\kappa - \frac{3\epsilon}{100}, \kappa+\frac{3\epsilon}{100}\right)} (F)$
	  (because both action windows contain the same set of periodic trajectories)
  \item 
	  $CF^{\left(\kappa - \frac{\epsilon}{100}, \kappa+\frac{\epsilon}{100}\right)} (F) \simeq 
	  HF^{\left(\kappa - \frac{\epsilon}{100}, \kappa+\frac{\epsilon}{100}\right)} (F)$
	  (since Floer differential vanishes inside the action filtration window)
  \item 
	  given a non-degenerate $f' \in \Ham (D)$ generated by a Hamiltonian $F'_t$ which is $\frac{\epsilon}{100}$-close to $F$, the continuation maps 
\[
  \begin{split}
	  \Psi:& CF^{\left(\kappa - \frac{\epsilon}{100}, \kappa+\frac{\epsilon}{100}\right)} (F) \to 
	  CF^{\left(\kappa - \frac{2\epsilon}{100}, \kappa+\frac{2\epsilon}{100}\right)} (F') \\
	  \widehat{\Psi}:& CF^{\left(\kappa - \frac{2\epsilon}{100}, \kappa+\frac{2\epsilon}{100}\right)} (F') \to 
	  CF^{\left(\kappa - \frac{3\epsilon}{100}, \kappa+\frac{3\epsilon}{100}\right)} (F) \simeq 
	  CF^{\left(\kappa - \frac{\epsilon}{100}, \kappa+\frac{\epsilon}{100}\right)} (F)
  \end{split}
\]
	  
	  compose to 
	  $$\widehat{\Psi} \circ \Psi: CF^{\left(\kappa - \frac{\epsilon}{100}, \kappa+\frac{\epsilon}{100}\right)} (F) \to 
	  CF^{\left(\kappa - \frac{\epsilon}{100}, \kappa+\frac{\epsilon}{100}\right)} (F)$$ 
	  which is chain homotopic to $\Id$. Moreover, the composition equals $\Id$ as Floer differential vanishes in the given action window.
\end{itemize}
From these facts the authors are able to derive stability of $B$ under non-degenerate perturbations with energy bounded by $\frac{\epsilon}{100}$
(the $\frac{1}{100}$ factor is not sharp).

In a similar manner, Definition~\ref{D:quasi-isolated} is interpreted in Proposition~3.4 as:
let $\kappa \in \mathcal{A}_F (B)$, denote by $B^\kappa \subset B$ the subset of trajectories whose action is $\kappa$. 
Then the vector space $\ZZ_2<B^\kappa>$ is a subcomplex of $CF^{\left(\kappa - \frac{\epsilon}{100}, \kappa+\frac{\epsilon}{100}\right)} (F)$
as orbits from $B$ are not involved in Floer trajectories with small energy. Given a non-degenerate Hamiltonian
$f' \in \Ham(\Sigma)$ which is generated by $F'_t$ and is $\frac{\epsilon}{100}$-close to $F$, the restricted continuation maps 
\[
  \begin{split}
	  \Psi:& \ZZ_2<B^\kappa> \to CF^{\left(\kappa - \frac{\epsilon}{100}, \kappa+\frac{\epsilon}{100}\right)} (F') \\
	  \widehat{\Psi}:& CF^{\left(\kappa - \frac{\epsilon}{100}, \kappa+\frac{\epsilon}{100}\right)} (F') \to \ZZ_2<B^\kappa>
  \end{split}
\]
  compose to $$\widehat{\Psi} \circ \Psi = \Id : \ZZ_2<B^\kappa> \to \ZZ_2<B^\kappa> .$$

In both cases the argument uses just the fact that the orbits of $B$ do not interact with other points in $\Fix{f}$ 
(and also do not interact between themselves) via `short' trajectories of Floer differential. Then a clever counting of intersections
for Floer cylinders of the continuation maps implies stability of the braid $B$ under the perturbation to $f'$. 
However, this technical condition of `non-interaction' in a narrow action window can be established by much less restrictive assumptions than
those in Definitions~\ref{D:isolated} or \ref{D:quasi-isolated}, leaving ample space for generalizations of the braid stability phenomenon.
Moreover, the geometric condition that the orbits of $B$ are not involved in `short' trajectories of the Floer differential
may be replaced with an algebraic one that the Floer differential vanishes within the action window of length $2 \epsilon$.
That is, short Floer trajectories may exist, but they show up in pairs and cancel out in the differential 
(we consider Floer theory with $\ZZ_2$ coefficients). This opens a way to extend the braid stability to certain scenarios with degenerate
$f \in \Ham(\Sigma)$, including many setups with autonomous and integrable diffeomorphisms. We explain below the (very minor)
modifications to the argument of ~\cite{Al-Me:braid} which are necessary to establish this generalization.

\begin{dfn} \label{D:admis}
  We call a braid $B$ formed by a possibly degenerate $f \in \Ham(D)$ \emph{$\epsilon$-admissible} if $B$ and $f$ satisfy the following conditions.
  Pick a compactly-supported Hamiltonian function $F_t$ generating $f$.
  \begin{itemize}
	  \item 
		  The fixed points of $f$ are either non-degenerate (in particular, appear as isolated points in $\Fix(f)$) or 
		  form isolated submanifolds in $D$, and the action is locally constant on each such submanifold.
	  \item
		  $\Spec(F) = \mathcal{A}_F (\Fix(f))$ is a finite set. For each $\kappa \in \Spec (F)$ there is just one connected component in $\Fix(f)$
		  (which may be a single non-degenerate point) whose action is $\kappa$. 
		  
		  This way, the only orbits with action zero are degenerate fixed points near the boundary.
	  \item
		  Different orbits in $B$ correspond to different connected components of fixed points (and therefore have different actions).
		  
		  Moreover, each $\gamma \in B$ corresponds either to an isolated fixed point or the appropriate submanifold $S_\gamma \subset \Fix(f)$ is 
		  diffeomorphic to $S^1$.
	  \item
		  $\Spec(F)$ is $\epsilon$-isolated. That is, for all distinct $\kappa_1, \kappa_2 \in \Spec(F)$ holds 
		  $|\kappa_1 - \kappa_2| \geq \epsilon$. In particular, fixed points near the boundary are separated in spectrum from the 
		  remaining $p \in \Fix(f)$.
   \end{itemize}
\end{dfn}

While these assumptions are rather restrictive and can be greatly relaxed, they are sufficient to establish the braid stability for the example
constructed in Section~\ref{S:Proof}.

\medskip 

\begin{rmk} \label{Rmk:pert}
Let $f \in \Ham(D)$ generated by $F_t$ and $B$ a braid of fixed points of $f$ which are $\epsilon$-admissible for some $\epsilon > 0$.
Pick $0 < \delta << \epsilon$.
Using standard perturbation techniques, $F_t$ can be perturbed in a neighborhood of $\Fix(f)$ to $\widehat{F_t}$ which generates a non-degenerate $\hat{f} \in \Ham(D)$ 
(non-degenerate in the interior of its support) such that
\begin{itemize} 
  \item 
	  the perturbation is $\delta$-small in both $C^\infty$ and Hofer metrics.
  \item
	  $\Fix (\hat{f}) \subset \Fix (f)$ and for each $p \in \Fix (\hat{f})$ the action difference 
	  $|\mathcal{A}_{\widehat{F}} (p) - \mathcal{A}_{F} (p)| < \delta$.
	\item
		Non-degenerate $p \in \Fix(f)$ are fixed points also for $\hat{f}$.
		If certain $\gamma \in B$ is induced by a degenerate fixed point $p \in \Fix(f)$, then $p \in \Fix(\hat{f})$. Moreover, the circle 
		$S_p \subset \Fix(f)$ of fixed points which contains $p$ splits after the perturbation to $\hat{f}$ into a pair of fixed points 
		$p, p' \in S_p \cap \Fix(\hat{f})$.
		The Conley-Zehnder indices of $p$ and $p'$ differ by one and for an appropriate choice of almost complex structure, 
		$p$ and $p'$ are connected in $CF(\widehat{F})$ by a pair of Floer trajectories of energy $< 2 \delta$.
\end{itemize}
\end{rmk}
This way if the braid $B$ was induced by orbits of fixed points $p_1, \ldots, p_k \in \Fix(f)$, after the perturbation
the same set of fixed points $p_1, \ldots, p_k \in \Fix(\hat{f})$ induces a braid $\widehat{B}$ whose strands are $\delta$-close to $B$
in the $C^0$ metric. Hence $B = \widehat{B}$ assuming $\delta$ is small enough. That is, $B$ `survives' the perturbation. 
We also remark that every braid generated by points in $\Fix(\hat{f})$ appears also in the original $f$ 
(and actually is induced by the same fixed points).

We sketch the proof for an analogue of Proposition~3.4 of ~\cite{Al-Me:braid} for the perturbed diffeomorphism $\hat{f}$.
\begin{prop} \label{P:continuation}
Let $f$, $B$ be $\epsilon$-admissible and $\hat{f} \in \Ham (D)$ be a non-degenerate $\delta$-close perturbation constructed
as described above. Suppose $f' \in \Ham(D)$ is non-degenerate and is generated by a Hamiltonian function
$F'_t$ which is $\epsilon'$-close to $F_t$ where $\epsilon' = \frac{\epsilon}{2} - \delta$. Let $\gamma \in B$ be an orbit generated by
$p \in \Fix (f)$.
Then the restricted continuation maps 
\[	
	\begin{split}
	  &\Psi : \ZZ_2<\gamma> \to CF^{\left(\kappa-\frac{\epsilon}{2}, \kappa+\frac{\epsilon}{2}\right)}(F') \\
	  &\widehat{\Psi} : CF^{\left(\kappa-\frac{\epsilon}{2}, \kappa+\frac{\epsilon}{2}\right)}(F') \to \ZZ_2<\gamma>\\
	\end{split}
\]
constructed as in Proposition~3.4, \cite{Al-Me:braid} satisfy
\[
	\widehat{\Psi} \circ \Psi = \Id : \ZZ_2<\gamma> \to \ZZ_2<\gamma>
\]
\end{prop}
Proof:
Suppose first that $p$ is degenerate for $f$. By Definition~\ref{D:admis}, it belongs to a circle $S_p \subset \Fix (f)$.
Then after the perturbation there are exactly two fixed points $p, p' \in \Fix(\hat{f})$ whose action falls within
$(\kappa - \delta, \kappa+\delta)$. Moreover, as $f, B$ are $\epsilon$-admissible, there are no other points from $\Fix(\hat{f})$ 
in the action interval $(\kappa - \epsilon + \delta, \kappa+\epsilon-\delta)$.

As the Floer differential in $CF^{(\kappa-\delta, \kappa+\delta)}(\widehat{F})$ vanishes
(there are just two Floer cylinders that connect the orbits $\gamma, \gamma_{p'}$ of $p$ and $p'$, and they cancel out), we have 
\[
	CF^{(\kappa-\delta, \kappa+\delta)}(\widehat{F}) \simeq CF^{(\kappa - \epsilon + \delta, \kappa+\epsilon-\delta)} \simeq \ZZ_2<\gamma, \gamma_{p'}> .
\]

We apply Proposition~3.8 from \cite{Al-Me:braid} and obtain
\[
  \begin{split}
	  \Psi:& CF^{(\kappa-\delta, \kappa+\delta)}(\widehat{F}) \to CF^{\left(\kappa-\frac{\epsilon}{2}, \kappa+\frac{\epsilon}{2}\right)}(F') \\
	  \widehat{\Psi}:& CF^{\left(\kappa-\frac{\epsilon}{2}, \kappa+\frac{\epsilon}{2}\right)}(F') \to 
	  CF^{(\kappa - \epsilon + \delta, \kappa+\epsilon-\delta)}(\widehat{F}) \simeq 
	  CF^{(\kappa-\delta, \kappa+\delta)}(\widehat{F})\\
  \end{split}
\]
such that the composition $\widehat{\Psi} \circ \Psi$ is chain homotopic to the identity map of $CF^{(\kappa-\delta, \kappa+\delta)}(\widehat{F})$.
Since the Floer differential vanishes, $\widehat{\Psi} \circ \Psi = \Id$.
The proposition follows by restricting the domain of $\Psi$ to $\ZZ_2<\gamma>$ and taking a quotient of the image of $\widehat{\Psi}$.

If $p$ is non-degenerate, the claim follows by a direct application of Proposition~3.8.
\Qed.

Now we are ready to prove
\begin{thm} \label{T:stability}
  Suppose a possibly degenerate $f \in \Ham(D)$ and a braid $B$ induced by points of $\Fix(f)$ satisfy the $\epsilon$-admissibility criterion.
  Then for all non-degenerate $f' \in \Ham(D)$ with $\dist_H(f, f') < \frac{\epsilon}{100}$, the fixed points of $f'$ generate a copy of $B$.
\end{thm}
Proof:
  Denote $\epsilon' = \dist_H(f, f')$ and pick generating functions $F_t, F'_t$ with $\dist_H(F_t, F'_t) < 2 \epsilon'$.
  Pick $\delta > 0$ such that $\delta + 2\epsilon' < \frac{\epsilon}{2}$ and construct a $\delta$-small perturbation 
  $\hat{f}$ generated by $\widehat{F}_t$ as in Remark~\ref{Rmk:pert}. As it was explained above, a copy of $B$ is induced also by the perturbed Hamiltonian $\hat{f}$.
  $\dist_H (\widehat{F}, F') < \delta + 2\epsilon' < \frac{\epsilon}{2} - \delta$, hence Proposition~\ref{P:continuation} holds.
  Now one proceeds as in the proof of Theorem~3 in ~\cite{Al-Me:braid}, using Proposition~\ref{P:continuation} 
  instead of Proposition~3.8 of \cite{Al-Me:braid} in Step~1 of the argument. The rest of the argument applies verbatim, showing that $B$ persists in $f'$.
\Qed

\begin{rmk}
	The bound $\frac{\epsilon}{100}$ is taken from ~\cite{Al-Me:braid} and is very far from being sharp. 
	The sharp estimate seems to be somewhere near $\frac{\epsilon}{2}$.
\end{rmk}

\section{Proof of Theorem~\ref{T:main} for the unit disk $D$}\label{S:Proof}
\begin{figure}[!htbp] 
\begin{center}
\includegraphics[width=0.9\textwidth]{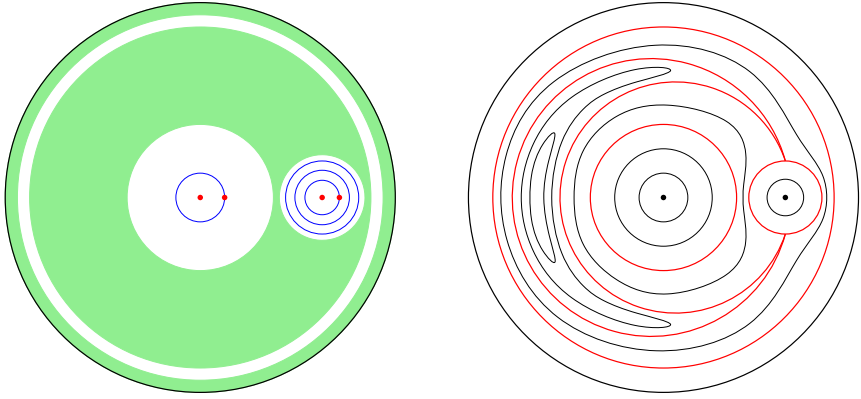}
\caption{The fixed points of $g$ and an invariant lamination.}  \label{F:fixed}
\end{center}
\end{figure}

Denote by $\Delta$ the disk of radius $0.3$ around $(0.5, 0)$ in $D$.
We build an integrable $g \in \Ham (D)$ as a composition of two autonomous rotations in $D$ and in $\Delta$.
Pick smooth non-increasing functions $\alpha : [0,  1] \to [0, 2.5]$ and $\beta: [0, 0.3] \to [3.5, 0]$
which satisfy:
\[
  \begin{split}
  &\left\{\begin{split}
	 &\alpha (0) > 2 \\
	 &\alpha (0.1) = 2 \\
	 &\alpha(r) = 1 \; \text{for} \; 0.2 \leq r \leq 0.8 \\
	 &\alpha (r) = 0 \; \text{for} \; r \geq 0.9 \\
	 &\alpha \; \text{is strongly decreasing in} \; [0, 0.2] \; \text{and} \; [0.8, 0.9],
  \end{split} \right. \\
  &\left\{\begin{split}
	 &\beta (0) > 3 \\
	 &\beta (0.1) = 3 \\
	 &\beta (r) = 0 \; \text{for} \; r \geq 0.2 \\
	 &\beta \; \text{is strongly decreasing in} \; [0, 0.2] .
  \end{split} \right.
  \end{split}
\]
\begin{figure}[!htbp]  
\begin{center}
\includegraphics[width=0.6\textwidth]{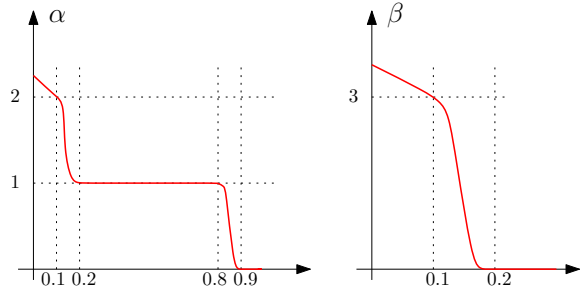}
\caption{$\alpha$ and $\beta$.}   \label{F:graph}
\end{center}
\end{figure}
Example graphs of such functions can be seen in Figure~\ref{F:graph}. 
Later on in the text we will impose additional restrictions on $\alpha, \beta$, which will be satisfied for a generic choice of functions.
Let $h \in \Ham (D)$ be the map which rotates concentric circles $\{x^2 + y^2 = r^2\} \subset D$ counterclockwise 
by the angle $2 \pi \cdot \alpha(r)$ and $h_\Delta$ be the map which rotates the circles $\{(x-0.5)^2 + y^2 = r^2\} \subset \Delta$ by
$2 \pi \cdot \beta(r)$. The map $h_\Delta$ is extended to the rest of $D$ by the identity.
It is easy to see that $h, h_{\Delta} \in \Ham_{Aut} (D)$ and they commute. Moreover, as they both preserve the lamination on
Figure~\ref{F:fixed}, the same is true for the composition $h_{\Delta} \circ h$. It is not difficult to construct an integral function $D \to \RR$ 
which defines this lamination by its level sets, hence $g := h_{\Delta} \circ h \in \Ham_{Int} (D)$.
\begin{figure}[!htbp] 
\begin{center}
\includegraphics[width=0.4\textwidth]{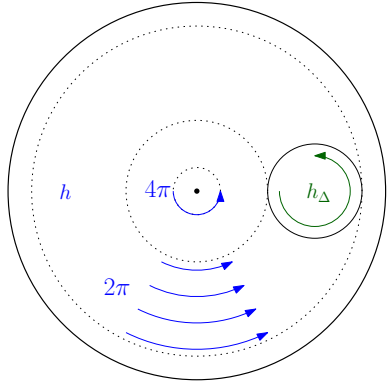}
\caption{$g = h_{\Delta} \circ h$.}  \label{F:diffeo}
\end{center}
\end{figure}

We claim that $g$ is not autonomous. Indeed, the points 
$$s = (0,0), \, p_1 = (0.1, 0), \, p_2 = (0.5, 0), \, m = (0, 0.6) \in \Fix (h) \cap \Fix (h_\Delta)$$
and by the choice of rotation functions $\alpha, \beta$
$$	\w_g (p_1, s) = \w_h (p_1, s) + \w_{h_\Delta} (p_1, s) = 2 + 0 = 2 . $$
In a similar way we compute
\[
  \begin{split}
	\w_g (p_2, s) &= \w_g (p_2, p_1) = 1, \\
	\w_g (m, s) &= \w_g (m, p_1) = 1, \quad \w_g (m, p_2) = 3 
 \end{split}
\]
-- the winding numbers are the same as for the strands of $B_{na}$. This implies that under any Hamiltonian isotopy 
$g_t$ connecting $g$ to $\Id$, the trajectories of $\{s, p_1, p_2, m\}$ satisfy the obstruction for autonomous braids. 
(With some extra efforts one may see that the braid $B$ they induce is indeed homotopic to $B_{na}$.)

\underline{Claim 1:} the braid $B$ is stable under Hofer-small perturbations of $g$ which are regular enough.

Indeed, we verify that $g$ and $B$ are $\epsilon$-admissible for a certain $\epsilon > 0$ (see Definition~\ref{D:admis}). 
The fixed points of $g$ are shown on Figure~\ref{F:fixed}: two two-dimensional families (green), one-dimensional families in blue 
and two non-degenerate fixed points. Red dots indicate $s, p_1, p_2, m$. 
A straightforward computation shows that the action is constant in each submanifold of degenerate fixed points.
By adjusting the functions $\alpha, \beta$ in their `flexible' domains one may ensure that the action values in each connected 
component are distinct and differ from the actions of non-degenerate fixed points. Therefore there exists an $\epsilon > 0$ such 
that $\Spec (g)$ is $\epsilon$-separated.
Applying Theorem~\ref{T:stability}, one sees that all non-degenerate $f \in \Ham (D)$ in the $\frac{\epsilon}{100}$-Hofer-neighborhood
of $g$ induce a copy of $B$.

\underline{Claim 2:} the $\frac{\epsilon}{100}$-Hofer neighborhood of $g$ does not intersect $\Ham_{Aut} (D)$.

Suppose by contradiction that there exists $h \in \Ham_{Aut} (D)$ such that $\dist_H (g, h) < \frac{\epsilon}{100}$. 
Denote by $H : D \to \RR$ a compactly supported time-independent function which generates $h$.
Applying a generic $C^\infty$-small compactly supported perturbation to $H$ we may assume that the resulting function $H'$ is 
Morse in the interior of its support and $1$-periodic orbits for the flow of $H'$ consist of the critical points of $H'$ and a finite collection of circles.
Assuming the perturbation is small enough, the Hamiltonian $h'$ generated by $H'$ is still $\frac{\epsilon}{100}$-close to $g$.
We perturb $h'$ farther to a non-degenerate $h'' \in \Ham (D)$ in a similar way to Remark~\ref{Rmk:pert}, forcing each circle of $1$-periodic points
of $h'$ to split into a pair of fixed points of $h''$, without introducing new periodic trajectories. As it is mentioned after Remark~\ref{Rmk:pert},
all braids induced by $h''$ appear also for $h'$. Given that the perturbation of $h'$ to $h''$ is sufficiently small, 
$\dist_H (g, h'') < \frac{\epsilon}{100}$ hence ($h''$ is non-degenerate) $h''$ induces a copy of $B$. But that means that $B$ appears
also for $h' \in \Ham_{Aut} (D)$, which is a contradiction as the braid $B$ is `non-autonomous'.

\section{Adaptation for general surfaces}\label{S:Surf}

\subsection{Other surfaces except $S^2$:} 
\subsubsection{Non-autonomous braids.}
Let $\Sigma$ be a compact symplectic surface, denote by $B^0_k(\Sigma)$ the subgroup of the pure braid group $B_k(\Sigma)$ 
formed by braids whose strands correspond to contractible loops. Pick $f \in \Ham (\Sigma)$.
$\Ham (\Sigma)$ is simply connected (see Chapter~7.2 in \cite{Po:geo-symp}), hence by Arnold conjecture each $p \in \Fix (f)$ defines a free homotopy class
$[f_t (p)]$ which does not depend on the isotopy $f_t$ that connects $f$ to the $\Id_\Sigma$. Denote by $\Fix_0 (f)$
the set of fixed points of $f$ with contractible trajectories.
Given a $k$-tuple of distinct points $\mathbf{p} = (p_1, \ldots, p_k) \in \Fix_0 (\Sigma)$ and
an isotopy $f_t \in \Ham (\Sigma)$ which connects $f$ to the identity, consider the braid $B_\mathbf{p} \subset B^0_k (\Sigma)$
defined by the trajectory $\{f_t (\mathbf{p})\}_t$.
As in the case of the disk, this definition is ambiguous, and this time the ambiguity might be more complicated than a mere conjugation.
The appropriate notion is a \emph{braid type} of $B_\mathbf{p}$ (see Chapter~2 in \cite{Mat:per-braid} or a short summary in Sections~1.3 and 1.4 of \cite{Al-Me:braid}).
Two braids $B, B' \in B^0_k (\Sigma)$ represent the same braid type if and only if they are freely isotopic as braids (\cite{Mat:per-braid}).
This allows to introduce topological invariants similar to the winding number which are well defined for a braid type.

Denote by $B_k^{0,Aut}(\Sigma) \subset B^0_k(\Sigma)$ the set of \emph{autonomous braids} that are induced by autonomous flows $g_t$.
We adapt the obstruction from Section~\ref{S:Braids} to general surfaces different from $S^2$. We note that the interior of the 
universal cover $\widetilde{\Sigma}$ is homeomorphic to a disk. Using this identification we introduce the following 
set-valued definition of the \emph{winding number}. Given a flow $g_t$ and two distinct points $p, q \in \Fix_0(g_1)$, let
\[
  \w_g (p, q) = \left\{\deg \left(\Arg (\tilde{g}_t(\tilde{p}) - \tilde{g}_t (\tilde{q}))\right) \right\}
\]
where $\tilde{g}_t$ is the lift of $g_t$ to a [possibly non-compactly supported] Hamiltonian isotopy on $\widetilde{\Sigma}$
which starts at $\Id_{\widetilde{\Sigma}}$ and $\tilde{p}, \tilde{q}$ run over all lifts of $p, q$ to $\widetilde{\Sigma}$.

Once again, the set-valued winding numbers descend to braids $B_\mathbf{p}$ induced by periodic trajectories of the flow.
They are not affected by the ambiguity in the definition of $B_\mathbf{p}$ as two braids of the same braid type are freely isotopic,
hence lifts of their strands to $\widetilde{\Sigma}$ are freely isotopic as well, thus have the same winding numbers.

\medskip

Suppose $g_t$ is an autonomous flow in $\Diff(\Sigma)$. 
We introduce a partial order on the set of contractible $1$-periodic trajectories of $g_t$:
$c \leq c'$ if $D_c \subset D_{c'}$ where $D_c, D_{c'}$ are the disks bounded by $c, c'$.
Or, in other words, there exists a lift of $c'$ which either encircles a lift of $c$ or coincides with it. 
(Once again, we agree that equilibrium points $q$ bound a singleton set $\{q\}$).
This defines a weak partial order.
One shows that $\w_g (p, q) = \{0\}$ if and only if the trajectories $c_p, c_q$ are incomparable. 
Similar to the case of $D$, if three distinct $1$-periodic points $p, q, q'$ satisfy $c_{q} \leq c_p$ and $c_{q'} \leq c_p$ then $\w_g (p, q) = \w (p, q')$.

Recall the example $B_{na} (D) \in B_4 (D)$ of a `non-autonomous' braid in $D$.
Using an embedding $D \hookrightarrow \Sigma$ we construct $B_{na} (\Sigma) \in B^0_4 (\Sigma)$. Assuming that $\Sigma$ is homeomorphic 
neither to $D$ nor to $S^2$, its fundamental group is not trivial and one shows
\[
 \begin{split}
	\w (p_1, s) &= \{2, 0\},  \\
	\w (p_2, s) &= \w (p_2, p_1) = \{1, 0\}, \\
	\w (m, s) &= \w (m, p_1) = \{1, 0\}, \quad \w (m, p_2) = \{3, 0\} .
 \end{split}
\]
As before, assuming by contradiction that $B_{na} (\Sigma)$ is induced by an autonomous flow,
this set of winding numbers implies that all trajectories are comparable under the partial order.
On the other hand, the winding numbers above obstruct existence of a maximal element, which gives a contradiction.
Hence $B_{na} (\Sigma)$ is `non-autonomous'.

\begin{rmk}
  There is even simpler example if one elects to work with braids of non-contractible loops:
  a `planet' traversing a simple non-contractible orbit and its `moon'. One can verify that it defines an integrable non-autonomous braid with 
  just two strands.
\end{rmk}

\subsubsection{Braid stability.}
\cite{Al-Me:braid} provides a separate argument for the braid stability in closed surfaces different from $S^2$. There are few technical differences 
from the case of $\Sigma = D$, mostly related to the need to accommodate braids of non-contractible loops. 
As long as we restrict our attention to braids induced by contractible trajectories, the argument is nearly the same as before. 
Moreover, using the same approach as for $\Sigma = D$, it extends also to compact surfaces with boundary and admits a generalization to
degenerate diffeomorphisms which are admissible. 

We briefly mention the main steps. 
In Definition~\ref{D:admis} of an \emph{$\epsilon$-admissible} pair we replace $\Fix (f)$ by $\Fix_0 (f)$.
In the case $\Sigma$ is closed we don't need anymore to separate fixed points near the boundary that have action zero from 
the remaining fixed points. The rest is exactly the same.
Perturbations defined in Remark~\ref{Rmk:pert} and Proposition~\ref{P:continuation} work exactly the same way.

\begin{thm} \label{T:sigma_stability}
  Suppose a possibly degenerate $f \in \Ham(\Sigma)$ and a braid $B$ induced by points of $\Fix_0(f)$ satisfy the $\epsilon$-admissibility criterion.
  Then for all non-degenerate $f' \in \Ham(\Sigma)$ with $\dist_H(f, f') < \frac{\epsilon}{100}$, the fixed points of $f'$ generate a copy of $B$.
\end{thm}
The proof repeats the proof of Theorem~\ref{T:stability}, invoking the argument from Section~6.2 of \cite{Al-Me:braid} instead of
Theorem~3 of the same paper.

\subsubsection{Construction.}
Without loss of generality, assume that $\Area (\Sigma) = \int_\Sigma \omega > \pi$ hence $D$ admits a symplectic embedding into $\Sigma$. 
In what follows we identify $D$ with its symplectic copy in $\Sigma$ and extend $g \in \Ham_{Int} (D)$
constructed in Section~\ref{S:Proof} to the rest of $\Sigma$ by the identity map. This extended map will be denoted by $g_\Sigma$.
Note that all fixed points of $g_\Sigma$ are contractible, and integrability of $g$ implies that the extended map is also integrable.
As $\{s, p_1, p_2, m\}$ defined a copy of $B_{na} (D)$ under $g$, they induce the non-autonomous braid $B_{na} (\Sigma)$ under $g_\Sigma$ 
(up to ambiguity of the braid type as discussed above).

Assuming the braid stability result for braids in $B^0_k (\Sigma)$, the rest of the argument is the same as in Section~\ref{S:Proof}, replacing
Theorem~\ref{T:stability} with Theorem~\ref{T:sigma_stability}.

\subsection{Proof for $S^2$:}
\subsubsection{Braids on $S^2$.}
Let $\Sigma = S^2$ be a sphere, denote by $B_k(S^2)$ the pure braid group of $S^2$. Pick $f \in \Ham (S^2)$ and 
an isotopy $f_t \in \Ham (S^2)$ connecting $f$ to $\Id_{S^2}$. 
Given a $k$-tuple of distinct points $\mathbf{p} = (p_1, \ldots, p_k) \in \Fix (f)$ consider the braid $B_\mathbf{p} \subset B_k (S^2)$
defined by the trajectory $\{f_t (\mathbf{p})\}$.
$\Ham (S^2)$ is not simply connected (its fundamental group is $\ZZ_2$, see Chapter~7.2 in \cite{Po:geo-symp}), hence in the construction of $B_\mathbf{p}$ from 
$f$ and $\mathbf{p}$ there is additional ambiguity derived from the choice of $\{f_t\}_t$. 

\subsubsection{Braid stability.}
In the stability argument the main difference is (compared to the previous cases) that the action of a trajectory in $S^2$ depends on a capping and in the 
Gromov compactness arguments one has to take into account the possibility of bubbling. Bubbling, however, cannot occur in trajectories 
of small energy. 
We assume that $\Area (S^2) = \int_{S^2} \omega = 1$ and will consider action values in $\RR / \ZZ$, in order to avoid dependence on the capping.

We briefly mention the main steps of the argument and the necessary modifications. 
In Definition~\ref{D:admis} of an \emph{$\epsilon$-admissible} pair we ask for $\Spec(f)$ to be finite and $\epsilon$-separated as a 
subset of $\RR / \ZZ$. There is no need to care about fixed points near the boundary. In addition we ask that $\epsilon << 1$ 
(in the text \cite{Al-Me:braid} $\epsilon < \frac{1}{3200}$, though this bound looks exaggerated.)
There are substantial differences in construction and behavior of various ingredients of the Floer theory, 
but once those are established, the proofs look very similar as long as energy constraints prevent bubbling. The
perturbations described in Remark~\ref{Rmk:pert} and the proof of Proposition~\ref{P:continuation} work exactly the same. 

\begin{thm} \label{T:s2_stability}
  Suppose a possibly degenerate $f \in \Ham(S^2)$ and a braid $B$ induced by points of $\Fix (f)$ satisfy the $\epsilon$-admissibility criterion.
  Then for all non-degenerate $f' \in \Ham(S^2)$ with $\dist_H(f, f') < \frac{\epsilon}{100}$, the fixed points of $f'$ generate a copy of $B$.
\end{thm}
The proof repeats the proof of Theorem~\ref{T:stability}, invoking the proof of Theorem~$1^*$ \cite{Al-Me:braid} instead of
Theorem~3 of the same paper.

\subsubsection{Construction.}
From the presentation above, the statement of Theorem~\ref{T:main} for spheres follows 
once one constructs a `non-autonomous' braid $B \in B_k (S^2)$ induced by an $\epsilon$-admissible $g \in \Ham_{Int} (S^2)$.
Indeed, assuming the braid stability result for $S^2$, the rest of the argument is the same as in Section~\ref{S:Proof}, replacing
Theorem~\ref{T:stability} with Theorem~\ref{T:s2_stability}.

Unfortunately, the author was not able to find an example of such a braid which is as simple as it was in the case of other surfaces. 
The proposed argument is based on a braid with 16 strands
which is definitely far from being the simplest representative in the class of `integrable non-autonomous' braids.
In order to avoid an inadequately long argument, we describe only a brief scheme of the proof.

\medskip

\underline{Step 1:} rescale the symplectic form on $D$ such that $\Area (D) = \int_D \omega = \frac{1}{10}$. 
Now one constructs four disjoint symplectic embeddings of $D$ into $S^2$. 
Denote by $D_1, D_2, D_3, D_4$ these four copies of $D$ in $S^2$ and construct $g_i \in \Ham_{Int} (D_i)$
as in Section~\ref{S:Proof}. These diffeomorphisms are extended in the complement of the four disks by the identity map. 
This extended map will be denoted by $g_{S^2}$. In the construction of $g_i$ we should use four different choices of 
the functions $\{\alpha_i, \beta_i\}$ to make sure that $\RR / \ZZ$-actions are different on each connected component of $\Fix(g_{S^2})$.
Integrability of each $g_i$ implies that the extended map $g_{S^2}$ is also integrable.
The four copies of $\{s, p_1, p_2, m\}$ in each $D_i$ induce under $g_{S^2}$ a braid $\check{B} \in B_{16}(S^2)$.
Since $\check{B}$ is clearly `integrable', it remains to show that it is `non-autonomous'.

\underline{Step 2:}
A $k$-tuple of loops $(\gamma_1, \ldots, \gamma_k)$ which define a braid on $\Sigma$ can be seen as a loop in the ordered configuration 
space $X_k (\Sigma)$. This results in identification of the pure braid group of $\Sigma$ with $\pi_1 (X_k (\Sigma))$.
In the case of $\Sigma = S^2$ and $k > 3$, $X_k (S^2)$ is homeomorphic to $X_3 (S^2) \times X_{k-3} (\CC \setminus \{0, 1\})$
(Theorem~2.1 in \cite{Fe-Zi:conf_spheres}). 

Using this decomposition, we send $\{(\gamma_1 (t), \ldots, \gamma_k (t))\}_t$ to the pair of braids 
$\{(\gamma_1 (t), \gamma_2 (t), \gamma_3 (t))\}_t$ and $\{(\gamma'_4 (t), \ldots, \gamma'_k (t))\}_t$
where $\gamma'_{3+i} (t)$ are obtained from $\gamma_{3+i} (t)$ by identifying $S^2 \setminus \{\gamma_1 (t), \gamma_2 (t), \gamma_3 (t)\}$
with $\CC \setminus \{0, 1\}$. That is, a selection of three strands in a braid $B \in B_k (S^2)$ induces a braid 
$B' \in B_{k-3} (\CC \setminus \{0, 1\})$.

\underline{Claim:}
given an `autonomous' $B \in B_k (S^2)$ with $k > 3$, there exist three strands of $B$ such that the induced braid 
$B' \in B_{k-3} (\CC \setminus \{0, 1\})$ is `autonomous' in $\CC \setminus \{0, 1\}$.

The proof is not hard but quite long. 
It uses a careful analysis of combinatorics of level sets of a function $H: S^2 \to \RR$.

\underline{Step 3:}
we claim that for any selection of three strands in $\check{B}$, the resulting braid $\check{B}' \in B_{13} (\CC \setminus \{0, 1\})$ 
is not autonomous. 
The idea is as follows. $\check{B}$ contains four independent copies of $B_{na} (D)$, removal of three strands will leave at least one copy 
undisrupted. One shows that it is mapped to a non-autonomous braid $B' \in B_4 (\CC \setminus \{0, 1\})$, which implies that $\check{B}'$ is not autonomous as 
well. (The proof follows from analysis of the winding (linking) numbers of strands in $B'$. One has to take into account that $B'$ might be 
defined by non-contractible loops, hence the obstruction used in the previous subsection has to be extended to cover this scenario.)

\bibliography{bibliography}

\end{document}